\documentclass[12pt]{article}
\usepackage{amsmath,amssymb,amsthm,latexsym,amscd}

\title{On semilattices and lattices \\ for families of
theories\footnote{
\newline\indent \ \ \ The research is partially supported by the Grants Council (under RF President) for State Aid of
Leading Scientific Schools (grant NSh-6848.2016.1) and by
Committee of Science in Education and Science Ministry of the
Republic of Kazakhstan (Grant No. 0830/GF4).} }
\author{Sergey V.
Sudoplatov\footnote{sudoplat@math.nsc.ru}}
\date{}
\begin{document}
\maketitle

\begin{abstract}
We define and study semilattices and lattices for $E$-closed
families of theories. Properties of these semilattices and
lattices are investigated. It is shown that lattices for families
of theories with least generating sets are distributive.

\end{abstract}

We continue to study structural properties of combin\-a\-tions of
structures and their theories \cite{cs, cl, lut, rest} defining
semilattices and lattices for families of theories. Properties of
these semilattices and lattices are investigated. It is shown that
lattices for families of theories with least generating sets are
distributive.

\section{Preliminaries}

Throughout the paper we use the following terminology in \cite{cs,
cl, rest}.

Let $P=(P_i)_{i\in I}$, be a family of nonempty unary predicates,
$(\mathcal{A}_i)_{i\in I}$ be a family of structures such that
$P_i$ is the universe of $\mathcal{A}_i$, $i\in I$, and the
symbols $P_i$ are disjoint with languages for the structures
$\mathcal{A}_j$, $j\in I$. The structure
$\mathcal{A}_P\rightleftharpoons\bigcup\limits_{i\in
I}\mathcal{A}_i$\index{$\mathcal{A}_P$} expanded by the predicates
$P_i$ is the {\em $P$-union}\index{$P$-union} of the structures
$\mathcal{A}_i$, and the operator mapping $(\mathcal{A}_i)_{i\in
I}$ to $\mathcal{A}_P$ is the {\em
$P$-operator}\index{$P$-operator}. The structure $\mathcal{A}_P$
is called the {\em $P$-combination}\index{$P$-combination} of the
structures $\mathcal{A}_i$ and denoted by ${\rm
Comb}_P(\mathcal{A}_i)_{i\in I}$\index{${\rm
Comb}_P(\mathcal{A}_i)_{i\in I}$} if
$\mathcal{A}_i=(\mathcal{A}_P\upharpoonright
A_i)\upharpoonright\Sigma(\mathcal{A}_i)$, $i\in I$. Structures
$\mathcal{A}'$, which are elementary equivalent to ${\rm
Comb}_P(\mathcal{A}_i)_{i\in I}$, will be also considered as
$P$-combinations.

Clearly, all structures $\mathcal{A}'\equiv {\rm
Comb}_P(\mathcal{A}_i)_{i\in I}$ are represented as unions of
their restrictions $\mathcal{A}'_i=(\mathcal{A}'\upharpoonright
P_i)\upharpoonright\Sigma(\mathcal{A}_i)$ if and only if the set
$p_\infty(x)=\{\neg P_i(x)\mid i\in I\}$ is inconsistent. If
$\mathcal{A}'\ne{\rm Comb}_P(\mathcal{A}'_i)_{i\in I}$, we write
$\mathcal{A}'={\rm Comb}_P(\mathcal{A}'_i)_{i\in
I\cup\{\infty\}}$, where
$\mathcal{A}'_\infty=\mathcal{A}'\upharpoonright
\bigcap\limits_{i\in I}\overline{P_i}$, maybe applying
Morleyzation. Moreover, we write ${\rm
Comb}_P(\mathcal{A}_i)_{i\in I\cup\{\infty\}}$\index{${\rm
Comb}_P(\mathcal{A}_i)_{i\in I\cup\{\infty\}}$} for ${\rm
Comb}_P(\mathcal{A}_i)_{i\in I}$ with the empty structure
$\mathcal{A}_\infty$.

Note that if all predicates $P_i$ are disjoint, a structure
$\mathcal{A}_P$ is a $P$-combination and a disjoint union of
structures $\mathcal{A}_i$. In this case the $P$-combination
$\mathcal{A}_P$ is called {\em
disjoint}.\index{$P$-combination!disjoint} Clearly, for any
disjoint $P$-combination $\mathcal{A}_P$, ${\rm
Th}(\mathcal{A}_P)={\rm Th}(\mathcal{A}'_P)$, where
$\mathcal{A}'_P$ is obtained from $\mathcal{A}_P$ replacing
$\mathcal{A}_i$ by pairwise disjoint
$\mathcal{A}'_i\equiv\mathcal{A}_i$, $i\in I$. Thus, in this case,
similar to structures the $P$-operator works for the theories
$T_i={\rm Th}(\mathcal{A}_i)$ producing the theory $T_P={\rm
Th}(\mathcal{A}_P)$\index{$T_P$}, being {\em
$P$-combination}\index{$P$-combination} of $T_i$, which is denoted
by ${\rm Comb}_P(T_i)_{i\in I}$.\index{${\rm Comb}_P(T_i)_{i\in
I}$}

For an equivalence relation $E$ replacing disjoint predicates
$P_i$ by $E$-classes we get the structure
$\mathcal{A}_E$\index{$\mathcal{A}_E$} being the {\em
$E$-union}\index{$E$-union} of the structures $\mathcal{A}_i$. In
this case the operator mapping $(\mathcal{A}_i)_{i\in I}$ to
$\mathcal{A}_E$ is the {\em $E$-operator}\index{$E$-operator}. The
structure $\mathcal{A}_E$ is also called the {\em
$E$-combination}\index{$E$-combination} of the structures
$\mathcal{A}_i$ and denoted by ${\rm Comb}_E(\mathcal{A}_i)_{i\in
I}$\index{${\rm Comb}_E(\mathcal{A}_i)_{i\in I}$}; here
$\mathcal{A}_i=(\mathcal{A}_E\upharpoonright
A_i)\upharpoonright\Sigma(\mathcal{A}_i)$, $i\in I$. Similar
above, structures $\mathcal{A}'$, which are elementary equivalent
to $\mathcal{A}_E$, are denoted by ${\rm
Comb}_E(\mathcal{A}'_j)_{j\in J}$, where $\mathcal{A}'_j$ are
restrictions of $\mathcal{A}'$ to its $E$-classes. The
$E$-operator works for the theories $T_i={\rm Th}(\mathcal{A}_i)$
producing the theory $T_E={\rm Th}(\mathcal{A}_E)$\index{$T_E$},
being {\em $E$-combination}\index{$E$-combination} of $T_i$, which
is denoted by ${\rm Comb}_E(T_i)_{i\in I}$\index{${\rm
Comb}_E(T_i)_{i\in I}$} or by ${\rm
Comb}_E(\mathcal{T})$\index{${\rm Comb}_E(\mathcal{T})$}, where
$\mathcal{T}=\{T_i\mid i\in I\}$.

Clearly, $\mathcal{A}'\equiv\mathcal{A}_P$ realizing $p_\infty(x)$
is not elementary embeddable into $\mathcal{A}_P$ and can not be
represented as a disjoint $P$-combination of
$\mathcal{A}'_i\equiv\mathcal{A}_i$, $i\in I$. At the same time,
there are $E$-combinations such that all
$\mathcal{A}'\equiv\mathcal{A}_E$ can be represented as
$E$-combinations of some $\mathcal{A}'_j\equiv\mathcal{A}_i$. We
call this representability of $\mathcal{A}'$ to be the {\em
$E$-representability}.

If there is $\mathcal{A}'\equiv\mathcal{A}_E$ which is not
$E$-representable, we have the $E'$-representability replacing $E$
by $E'$ such that $E'$ is obtained from $E$ adding equivalence
classes with models for all theories $T$, where $T$ is a theory of
a restriction $\mathcal{B}$ of a structure
$\mathcal{A}'\equiv\mathcal{A}_E$ to some $E$-class and
$\mathcal{B}$ is not elementary equivalent to the structures
$\mathcal{A}_i$. The resulting structure $\mathcal{A}_{E'}$ (with
the $E'$-representability) is a {\em
$e$-completion}\index{$e$-completion}, or a {\em
$e$-saturation}\index{$e$-saturation}, of $\mathcal{A}_{E}$. The
structure $\mathcal{A}_{E'}$ itself is called {\em
$e$-complete}\index{Structure!$e$-complete}, or {\em
$e$-saturated}\index{Structure!$e$-saturated}, or {\em
$e$-universal}\index{Structure!$e$-universal}, or {\em
$e$-largest}\index{Structure!$e$-largest}.

For a structure $\mathcal{A}_E$ the number of {\em
new}\index{Structure!new} structures with respect to the
structures $\mathcal{A}_i$, i.~e., of the structures $\mathcal{B}$
which are pairwise elementary non-equivalent and elementary
non-equivalent to the structures $\mathcal{A}_i$, is called the
{\em $e$-spectrum}\index{$e$-spectrum} of $\mathcal{A}_E$ and
denoted by $e$-${\rm Sp}(\mathcal{A}_E)$.\index{$e$-${\rm
Sp}(\mathcal{A}_E)$} The value ${\rm sup}\{e$-${\rm
Sp}(\mathcal{A}'))\mid\mathcal{A}'\equiv\mathcal{A}_E\}$ is called
the {\em $e$-spectrum}\index{$e$-spectrum} of the theory ${\rm
Th}(\mathcal{A}_E)$ and denoted by $e$-${\rm Sp}({\rm
Th}(\mathcal{A}_E))$.\index{$e$-${\rm Sp}({\rm
Th}(\mathcal{A}_E))$}

If $\mathcal{A}_E$ does not have $E$-classes $\mathcal{A}_i$,
which can be removed, with all $E$-classes
$\mathcal{A}_j\equiv\mathcal{A}_i$, preserving the theory ${\rm
Th}(\mathcal{A}_E)$, then $\mathcal{A}_E$ is called {\em
$e$-prime}\index{Structure!$e$-prime}, or {\em
$e$-minimal}\index{Structure!$e$-minimal}.

For a structure $\mathcal{A}'\equiv\mathcal{A}_E$ we denote by
${\rm TH}(\mathcal{A}')$ the set of all theories ${\rm
Th}(\mathcal{A}_i)$\index{${\rm Th}(\mathcal{A}_i)$} of
$E$-classes $\mathcal{A}_i$ in $\mathcal{A}'$.

By the definition, an $e$-minimal structure $\mathcal{A}'$
consists of $E$-classes with a minimal set ${\rm
TH}(\mathcal{A}')$. If ${\rm TH}(\mathcal{A}')$ is the least for
models of ${\rm Th}(\mathcal{A}')$ then $\mathcal{A}'$ is called
{\em $e$-least}.\index{Structure!$e$-least}

\medskip
{\bf Definition} \cite{cl}. Let $\overline{\mathcal{T}}$ be the
class of all complete elementary theories of relational languages.
For a set $\mathcal{T}\subset\overline{\mathcal{T}}$ we denote by
${\rm Cl}_E(\mathcal{T})$ the set of all theories ${\rm
Th}(\mathcal{A})$, where $\mathcal{A}$ is a structure of some
$E$-class in $\mathcal{A}'\equiv\mathcal{A}_E$,
$\mathcal{A}_E={\rm Comb}_E(\mathcal{A}_i)_{i\in I}$, ${\rm
Th}(\mathcal{A}_i)\in\mathcal{T}$. As usual, if $\mathcal{T}={\rm
Cl}_E(\mathcal{T})$ then $\mathcal{T}$ is said to be {\em
$E$-closed}.\index{Set!$E$-closed}

The operator ${\rm Cl}_E$ of $E$-closure can be naturally extended
to the classes $\mathcal{T}\subset\overline{\mathcal{T}}$ as
follows: ${\rm Cl}_E(\mathcal{T})$ is the union of all ${\rm
Cl}_E(\mathcal{T}_0)$ for subsets
$\mathcal{T}_0\subseteq\mathcal{T}$.

For a set $\mathcal{T}\subset\overline{\mathcal{T}}$ of theories
in a language $\Sigma$ and for a sentence $\varphi$ with
$\Sigma(\varphi)\subseteq\Sigma$ we denote by
$\mathcal{T}_\varphi$\index{$\mathcal{T}_\varphi$} the set
$\{T\in\mathcal{T}\mid\varphi\in T\}$.

\medskip
{\bf Proposition 1.1} \cite{cl}. {\em If
$\mathcal{T}\subset\overline{\mathcal{T}}$ is an infinite set and
$T\in\overline{\mathcal{T}}\setminus\mathcal{T}$ then $T\in{\rm
Cl}_E(\mathcal{T})$ {\rm (}i.e., $T$ is an {\sl accumulation
point} for $\mathcal{T}$ with respect to $E$-closure ${\rm
Cl}_E${\rm )} if and only if for any formula $\varphi\in T$ the
set $\mathcal{T}_\varphi$ is infinite.}

\medskip
{\bf Theorem 1.2} \cite{cl}. {\em For any sets
$\mathcal{T}_0,\mathcal{T}_1\subset\overline{\mathcal{T}}$, ${\rm
Cl}_E(\mathcal{T}_0\cup\mathcal{T}_1)={\rm
Cl}_E(\mathcal{T}_0)\cup{\rm Cl}_E(\mathcal{T}_1)$.}

\medskip
{\bf Definition} \cite{cl}. Let $\mathcal{T}_0$ be a closed set in
a topological space $(\mathcal{T},\mathcal{O}_E(\mathcal{T}))$,
where $\mathcal{O}_E(\mathcal{T})=\{\mathcal{T}\setminus{\rm
Cl}_E(\mathcal{T}')\mid\mathcal{T}'\subseteq\mathcal{T}\}$. A
subset $\mathcal{T}'_0\subseteq\mathcal{T}_0$ is said to be {\em
generating}\index{Set!generating} if $\mathcal{T}_0={\rm
Cl}_E(\mathcal{T}'_0)$. The generating set $\mathcal{T}'_0$ (for
$\mathcal{T}_0$) is {\em minimal}\index{Set!generating!minimal} if
$\mathcal{T}'_0$ does not contain proper generating subsets. A
minimal generating set $\mathcal{T}'_0$ is {\em
least}\index{Set!generating!least} if $\mathcal{T}'_0$ is
contained in each generating set for $\mathcal{T}_0$.

\medskip
{\bf Theorem 1.3} \cite{cl}. {\em If $\mathcal{T}'_0$ is a
generating set for a $E$-closed set $\mathcal{T}_0$ then the
following conditions are equivalent:

$(1)$ $\mathcal{T}'_0$ is the least generating set for
$\mathcal{T}_0$;

$(2)$ $\mathcal{T}'_0$ is a minimal generating set for
$\mathcal{T}_0$;

$(3)$ any theory in $\mathcal{T}'_0$ is isolated by some set
$(\mathcal{T}'_0)_\varphi$, i.e., for any $T\in\mathcal{T}'_0$
there is $\varphi\in T$ such that
$(\mathcal{T}'_0)_\varphi=\{T\}$;

$(4)$ any theory in $\mathcal{T}'_0$ is isolated by some set
$(\mathcal{T}_0)_\varphi$, i.e., for any $T\in\mathcal{T}'_0$
there is $\varphi\in T$ such that
$(\mathcal{T}_0)_\varphi=\{T\}$.}

\medskip
{\bf Definition} \cite{rest}. Two theories $T_1$ and $T_2$ of a
language $\Sigma$ are {\em disjoint}\index{Theories!disjoint
modulo $\Sigma_0$} modulo $\Sigma_0$, where
$\Sigma_0\subseteq\Sigma$, or {\em
$\Sigma_0$-disjoint}\index{Theories!$\Sigma_0$-disjoint} if $T_1$
and $T_2$ are do not have common nonempty predicates for
$\Sigma\setminus\Sigma_0$. If $T_1$ and $T_2$ are
$\varnothing$-disjoint, these theories are called simply {\em
disjoint}\index{Theories!disjoint}.

\section{Semilattices and lattices for families of theories}

{\bf Definition.} Let $X$ be a nonempty set of $E$-closed families
$\mathcal{T}\subset\overline{\mathcal{T}}$. Operations
$\mathcal{T}_1\wedge\mathcal{T}_2\rightleftharpoons\mathcal{T}_1\cap\mathcal{T}_2$
and $\mathcal{T}_1\vee\mathcal{T}_2\rightleftharpoons{\rm
Cl}_E(\mathcal{T}_1\cup\mathcal{T}_2)$, for $E$-closed
$\mathcal{T}_1,\mathcal{T}_2\subset\overline{\mathcal{T}}$,
generate a set $Y$ and form the structure $\langle Y;\,\wedge,
\vee\rangle$ denoted by $L(X)$.\index{$L(X)$}

\medskip
It is well known \cite{BurSan} that any $L(X)$ is a lattice
extensible to a complete lattice ${\rm CL}(X)$\index{${\rm
CL}(X)$} with
$$\bigwedge\limits_{j\in J}{\rm
Cl}_E(\mathcal{T}_j)=\bigcap\limits_{j\in J}{\rm
Cl}_E(\mathcal{T}_j)$$ and
$$
\bigvee\limits_{j\in J}{\rm Cl}_E(\mathcal{T}_j)={\rm
Cl}_E\left(\bigcup\limits_{j\in J}\mathcal{T}_j\right).
$$

By Theorem 1.2, for $E$-closed $\mathcal{T}_1,\mathcal{T}_2$,
${\rm Cl}_E(\mathcal{T}_0\cup\mathcal{T}_1)=\mathcal{T}_0\cup
\mathcal{T}_1$, i.e., the operation $\vee$ is the set-theoretic
union. At the same time, in general case, for $E$-closed
$\mathcal{T}_j$, $\bigvee\limits_{j\in
J}\mathcal{T}_j\ne\bigcup\limits_{j\in J}\mathcal{T}_j$, since,
for instance, the union of infinite set of singletons can generate
new theories. Thus, $L(X)$ is just a standard algebra with usual
set-theoretic unions and intersections (but can be without even
relative complements since these complements can be not
$E$-closed), whereas ${\rm CL}(X)$ is its natural
expansion-extension.

Now we consider restrictions of $L(X)$ in the following way.

For a nonempty set $X$ of $E$-closed families with least
generating sets, the operation $\vee$ generates a set $Z\subseteq
Y$ and forms a upper semilattice ${\rm SLLGS}(X)=\langle Z;\,\vee
\rangle$\index{${\rm SLLGS}(X)$} restricting the universe and the
language of $L(X)$.

Below we will show that ${\rm SLLGS}(X)$ always consists of
families with least generating sets whereas the operation $\wedge$
can generate a family without least generating sets.

\medskip
{\bf Proposition 2.1.} {\em If $E$-closed sets $\mathcal{T}_1$ and
$\mathcal{T}_2$, in a language $\Sigma$, have least generating
sets, then $\mathcal{T}_1\cup\mathcal{T}_2$, being $E$-closed, has
the least generating set.}

\medskip
{\bf\em Proof.} Let $\mathcal{T}'_1$ and $\mathcal{T}'_2$ be least
generating sets for $\mathcal{T}_1$ and $\mathcal{T}_2$
respectively, and $\mathcal{T}'_0$ be a subset of
$\mathcal{T}'_1\cup\mathcal{T}'_2$ consisting of all isolated
points with respect to $\mathcal{T}'_1\cup\mathcal{T}'_2$, i.e.,
of elements $T\in\mathcal{T}'_1\cup\mathcal{T}'_2$ with formulas
$\varphi\in T$ such that
$(\mathcal{T}'_1\cup\mathcal{T}'_2)_\varphi$ is a singleton.

Now we assume on contrary that $\mathcal{T}_1\cup\mathcal{T}_2$
does not have the least generating set. Then there is a theory
$T_1\in \mathcal{T}'_1\cup\mathcal{T}'_2$ such that $T_1\notin{\rm
Cl}_E(\mathcal{T}'_0)$. Without loss of generality we assume that
$T_1\in\mathcal{T}'_1$. Since $T_1$ is isolated with respect to
$\mathcal{T}'_1$ and not isolated with respect to
$\mathcal{T}'_1\cup\mathcal{T}'_2$ there is a formula $\varphi\in
T_1$ such that for any $\psi\in T_1$ forcing $\varphi$,
$(\mathcal{T}'_1)_\psi=\{T_1\}$ and
$(\mathcal{T}'_1\cup\mathcal{T}'_2)_\psi$ is infinite. Since
$T_1\notin{\rm Cl}_E(\mathcal{T}'_0)$, there are infinitely many
theories $T\in(\mathcal{T}'_2)_\psi$ which are not isolated with
respect to $\mathcal{T}'_1\cup\mathcal{T}'_2$. It implies that for
any formula $\chi\in T$ forcing $\psi$ there are infinitely many
theories in $(\mathcal{T}'_1)_\chi$. But since $\chi\vdash\psi$,
$(\mathcal{T}'_1)_\psi$ is infinite contradicting
$|(\mathcal{T}'_1)_\psi|=1$.~$\Box$

\medskip
{\bf Remark 2.2.} Arguments for \cite[Proposition 3.9]{rest} show
that the converse for Proposition 2.1 is not true, since there is
$\mathcal{T}_1\cup\mathcal{T}_2$ with the least generating set
such that $\mathcal{T}_1$ has the least generating set (for
$\mathcal{F}_q$ in terms of \cite{lut}) and $\mathcal{T}_2$ does
not have the least generating set (for $\{J_q\mid q\in\mathbb Q\}$
in terms of \cite{lut}).

If we denote by $\Sigma_0$ the set of nonempty predicates for
$\{J_q\mid q\in\mathbb Q\}$ and take a $\Sigma_0$-disjoint copy
$\mathcal{F}'_q$ for $\mathcal{F}_q$, which also generates
$\{J_q\mid q\in\mathbb Q\}$ with $J_q=\underline{\rm
lim}\,F_q=\underline{\rm lim}\,F'_q$, we get families
$\mathcal{T}$ and $\mathcal{T}'$ for $\{J_q\mid q\in\mathbb
Q\}\cup \mathcal{F}_q$ and $\{J_q\mid q\in\mathbb
Q\}\cup\mathcal{F}'_q$ respectively such that
$\mathcal{T}\cap\mathcal{T}'$ is a family of theories for
$\{J_q\mid q\in\mathbb Q\}$, which does not have the least
generating set.

\medskip
{\bf Remark 2.3.} The infinite semilattices ${\rm SLLGS}(X)$ can
be both complete and incomplete, and in the incomplete case ${\rm
SLLGS}(X)$ can not be extended to a complete semilattice
consisting of families with least generating sets.

Indeed, taking infinitely many $\Sigma_0$-disjoint copies
$\mathcal{F}^\mu_q$ of $\mathcal{F}_q$ \cite{lut}, $\mu<\lambda$,
and forming the set $X$ by $E$-closed families of theories for
$\{J_q\mid q\in\mathbb Q\}\cup\mathcal{F}^\mu_q$ \cite{lut} we can
freely unite elements of $X$ obtaining $E$-closed families with
least generating sets corresponding to unions of
$\mathcal{F}^\mu_q$.

At the same time, each singleton $\{T\}$, for
$T\in\overline{\mathcal{T}}$ is $E$-closed and with the least
generating set. Then taking a set $X$ of singletons we generate
the semilattice ${\rm SLLGS}(X)$ (which is in fact a distributive
lattice with related complements) consisting of all finite subsets
of $\cup X$. As there are $E$-closed families $\mathcal{T}$
without least generating sets, taking a (infinite) union of
singletons $\{T\}$ for $T\in\mathcal{T}$ we form the family
$\mathcal{T}$. Thus, infinite unions of families with least
generating sets can be without least generating sets, and in this
case ${\rm SLLGS}(X)$ can not be extended to a complete
semilattice consisting of families with least generating sets.

\medskip
Summarizing Proposition 2.1 and Remarks 2.2, 2.3 we have

\medskip
{\bf Theorem 2.3.} {\em $1.$ For any nonempty set $X$ of
$E$-closed families with least generating sets the structure ${\rm
SLLGS}(X)$ is a upper semilattice.

$2.$ There is a upper semilattice ${\rm SLLGS}(X)$ with elements
$x_1,x_2\in X$ having least generating sets and such that $x_1\cap
x_2$ does not have the least generating set.

$3.$ There is a upper semilattice ${\rm SLLGS}(X)$ which can not
be extended to a complete semilattice consisting of families with
least generating sets.}

\medskip
Now we take a nonempty set $X$ of $E$-closed families with least
generating sets and $\mathcal{T}_1, \mathcal{T}_2\in X$ with least
generating sets $\mathcal{T}'_1$ and $\mathcal{T}'_2$
respectively. We denote by $\mathcal{T}_1\wedge'
\mathcal{T}_2$\index{$\mathcal{T}_1\wedge' \mathcal{T}_2$} the
family $\mathcal{T}_0\subseteq \mathcal{T}_1\cap \mathcal{T}_2$
with the greatest generating set $\mathcal{T}'_0$ consisting of
all isolated points for $\mathcal{T}_1\cap \mathcal{T}_2$.

\medskip
{\bf Remark 2.4.} By the definition,
$\mathcal{T}'_0\supseteq(\mathcal{T}'_1\cap
\mathcal{T}_2)\cup(\mathcal{T}_1\cap \mathcal{T}'_2)$. In
particular, $\mathcal{T}'_0\supseteq\mathcal{T}'_1\cap
\mathcal{T}'_2$. These inclusions can be strict.

Indeed, take $\Sigma_0$-disjoint families $\mathcal{T}_1$ and
$\mathcal{T}_2$ of theories for $F_q\cup\{J_q\}$ and
$F'_q\cup\{J_q\}$ \cite{lut}, where $J_q=\underline{\rm
lim}\,F_q=\underline{\rm lim}\,F'_q$ and $\Sigma_0$ is the set of
predicate symbols interpreted by nonempty relations for $J_q$. For
the theory $T_0$ corresponding to $J_q$, we have
$\{T_0\}=\mathcal{T}_1\wedge' \mathcal{T}_2$, whereas
$\mathcal{T}'_1\cap \mathcal{T}'_2=\varnothing$.

\medskip
For the set $X$ the operations $\wedge'$ and $\vee$ generate a set
$U\supseteq X$ with a structure ${\rm
LLGS}(X)\rightleftharpoons\langle U;\,\wedge',\vee\rangle$.

Directly checking we have

\medskip
{\bf Proposition 2.5.} {\em Any structure ${\rm LLGS}(X)$ is a
lattice.}

\medskip

By the definition for every $\mathcal{T}_1, \mathcal{T}_2\in U$
with least generating sets $\mathcal{T}'_1$ and $\mathcal{T}'_2$
respectively, we have, in ${\rm LLGS}(X)$, that $\mathcal{T}_1\leq
\mathcal{T}_2$ if and only if $\mathcal{T}'_2$ consists of three
disjoint parts $\mathcal{T}'_{2,1}$, $\mathcal{T}'_{2,2}$,
$\mathcal{T}'_{2,3}$ such that:

1) $\mathcal{T}'_{2,1}\subseteq \mathcal{T}'_1$,

2) $(\mathcal{T}'_{2,2}\cup\mathcal{T}'_{2,3})\cap
\mathcal{T}'_1=\varnothing$,

3) $\mathcal{T}'_{2,2}$ is used for generations of elements in
$\mathcal{T}'_1\setminus\mathcal{T}'_2$;

4) $\mathcal{T}'_{2,3}$ is not used for generations of elements in
$\mathcal{T}'_1\setminus\mathcal{T}'_2$.

The following proposition is obvious.

\medskip
{\bf Proposition 2.6.} {\em If $\mathcal{T}'_2\setminus
\mathcal{T}'_1$ is finite then $\mathcal{T}'_2=\mathcal{T}'_1\cup
\mathcal{T}'_{2,3}$ and, moreover,
$\mathcal{T}_2=\mathcal{T}_1\cup \mathcal{T}'_{2,3}$.}

\medskip
{\bf Remark 2.7.} The set $\mathcal{T}'_{2,3}$ can have an
arbitrary cardinality whereas for each element in
$\mathcal{T}'_1\setminus\mathcal{T}'_2$ being isolated by some
$(\mathcal{T}'_1)_\varphi$, the neighbourhood
$(\mathcal{T}'_{2,2})_\varphi$ should contain infinitely many
elements and the cardinality of $\mathcal{T}'_{2,2}$ is at least
$\lambda\cdot\omega$, where $\lambda$ is the number of disjoint
$\mathcal{T}'_{2,2}$-neighbourhoods for the elements in
$\mathcal{T}'_1\setminus\mathcal{T}'_2$.

\medskip
{\bf Theorem 2.8.} {\em Any lattice ${\rm LLGS}(X)$ is
distributive.}

\medskip
{\bf\em Proof.} We have to show two identities:
\begin{equation}\label{eqrest5}
\mathcal{T}_1\wedge'(\mathcal{T}_2\cup\mathcal{T}_3)=(\mathcal{T}_1\wedge'\mathcal{T}_2)\cup(\mathcal{T}_1\wedge'\mathcal{T}_3),
\end{equation}
\begin{equation}\label{eqrest6}
\mathcal{T}_1\cup(\mathcal{T}_2\wedge'\mathcal{T}_3)=(\mathcal{T}_1\cup\mathcal{T}_2)\wedge'(\mathcal{T}_1\cup\mathcal{T}_3)
\end{equation}
for any $\mathcal{T}_1,\mathcal{T}_2,\mathcal{T}_3\in U$ with
least generating sets
$\mathcal{T}'_1,\mathcal{T}'_2,\mathcal{T}'_3$ respectively.

For the proof of (\ref{eqrest5}) we note that the least generating
set for $\mathcal{T}_1\wedge'(\mathcal{T}_2\cup\mathcal{T}_3)$
consists of isolated points $T$ belonging both to $\mathcal{T}_1$
and to $\mathcal{T}_2\cup\mathcal{T}_3$. But then $T$ belongs to
$\mathcal{T}_2$ or to $\mathcal{T}_3$. In the first case, $T$ is
an isolated point for $\mathcal{T}_1\wedge'\mathcal{T}_2$ and, in
the second case,~--- an isolated point for
$\mathcal{T}_1\wedge'\mathcal{T}_3$. Therefore,
$T\in(\mathcal{T}_1\wedge'\mathcal{T}_2)\cup(\mathcal{T}_1\wedge'\mathcal{T}_3)$
and thus
$\mathcal{T}_1\wedge'(\mathcal{T}_2\cup\mathcal{T}_3)\subseteq(\mathcal{T}_1\wedge'\mathcal{T}_2)\cup(\mathcal{T}_1\wedge'\mathcal{T}_3)$.

Conversely, if an isolated point $T$ belongs to
$(\mathcal{T}_1\wedge'\mathcal{T}_2)\cup(\mathcal{T}_1\wedge'\mathcal{T}_3)$
then $T$ belongs to $\mathcal{T}_1\wedge'\mathcal{T}_2$ or to
$\mathcal{T}_1\wedge'\mathcal{T}_3$. If $T\in
\mathcal{T}_1\wedge'\mathcal{T}_2$ then either $T$ is an isolated
point for $\mathcal{T}_1\wedge'\mathcal{T}_2$ or belong to the
$E$-closure of isolated points in
$\mathcal{T}_1\wedge'\mathcal{T}_2$. Anyway,
$T\in\mathcal{T}_1\wedge'(\mathcal{T}_2\cup\mathcal{T}_3)$.
Similarly we get
$T\in\mathcal{T}_1\wedge'(\mathcal{T}_2\cup\mathcal{T}_3)$ for any
isolated $T\in\mathcal{T}_1\wedge'\mathcal{T}_2$. Thus,
$(\mathcal{T}_1\wedge'\mathcal{T}_2)\cup(\mathcal{T}_1\wedge'\mathcal{T}_3)\subseteq\mathcal{T}_1\wedge'(\mathcal{T}_2\cup\mathcal{T}_3)$
and the identity (\ref{eqrest5}) holds.

For the proof of (\ref{eqrest6}) we note that the least generating
set for $\mathcal{T}_1\cup(\mathcal{T}_2\wedge'\mathcal{T}_3)$
consists of isolated points $T$ belonging to $\mathcal{T}_1$ and
being an isolated point for $\mathcal{T}_1$, or belonging to to
$\mathcal{T}_2\wedge'\mathcal{T}_3$, being an isolated point for
$\mathcal{T}_2\wedge'\mathcal{T}_3$, and then belonging to
$\mathcal{T}_2$ and to $\mathcal{T}_3$. If $T\in\mathcal{T}_1$
then $T\in\mathcal{T}_1\cup\mathcal{T}_2$ and
$T\in\mathcal{T}_1\cup\mathcal{T}_3$, whence $T\in
(\mathcal{T}_1\cup\mathcal{T}_2)\wedge'(\mathcal{T}_1\cup\mathcal{T}_3)$.
If $T\in\mathcal{T}_2\wedge'\mathcal{T}_3$ then again
$T\in\mathcal{T}_1\cup\mathcal{T}_2$ and
$T\in\mathcal{T}_1\cup\mathcal{T}_3$ implying $T\in
(\mathcal{T}_1\cup\mathcal{T}_2)\wedge'(\mathcal{T}_1\cup\mathcal{T}_3)$.

Conversely, if an isolated point $T$ belongs to
$(\mathcal{T}_1\cup\mathcal{T}_2)\wedge'(\mathcal{T}_1\cup\mathcal{T}_3)$
then $T$ belongs to $\mathcal{T}_1\cup\mathcal{T}_2$ and to
$\mathcal{T}_1\cup\mathcal{T}_3$. So $T\in \mathcal{T}_1$ or
$T\in\mathcal{T}_2$, and $T\in \mathcal{T}_1$ or
$T\in\mathcal{T}_2$. If $T\in\mathcal{T}_1$ then
$T\in\mathcal{T}_1\cup(\mathcal{T}_2\wedge'\mathcal{T}_3)$.
Otherwise, $T\in\mathcal{T}_2\wedge'\mathcal{T}_3$ and so again
$T\in\mathcal{T}_1\cup(\mathcal{T}_2\wedge'\mathcal{T}_3)$.

Thus, the identity (\ref{eqrest6}) holds. $\Box$

\medskip
{\bf Remark 2.9.} For every $E$-closed family $\mathcal{T}$ with
the least generating set $\mathcal{T}'$ there is a superatomic
Boolean algebra
$\mathcal{B}(\mathcal{T})$\index{$\mathcal{B}(\mathcal{T})$}
\cite{Day} consisting of all subsets of $\mathcal{T}$ generated by
arbitrary subsets of $\mathcal{T}'$. If $\mathcal{T}_1\leq
\mathcal{T}_2$ in $\mathcal{B}(\mathcal{T})$ and $\mathcal{T}_1,
\mathcal{T}_2$ have the least generating sets $\mathcal{T}'_1$ and
$\mathcal{T}'_2$, respectively, then
$\mathcal{T}'_1\subseteq\mathcal{T}'_2$ and vice versa. Thus,
$\mathcal{B}(\mathcal{T})$ is isomorphic to the Boolean algebra
$\mathcal{B}(\mathcal{T}')$ with the natural relation $\subseteq$
and consisting of all subsets of $\mathcal{T}'$.

\end{document}